\documentclass{article}

\newtheorem{fed}{\textbf{Definition}}[section]
\newtheorem{theorem}[fed]{\textbf{Theorem}}
\newtheorem{lemma}[fed]{\textbf{Lemma}}

\newtheorem{cor}[fed]{\textbf{Corollary}}
\newtheorem{conj}[fed]{\textbf{Conjecture}}
\usepackage{amssymb,bbm,graphicx,epsfig,psfrag,epic,eepic,latexsym}
\usepackage{amsmath}
\usepackage{mathrsfs}
\usepackage[dvips]{color}

\begin{document}
\title{Smooth nondisplaceability for fixed point sets of involutions}
\author{Urs Frauenfelder}
\maketitle
\begin{abstract}
We prove that on closed manifolds of odd Euler characteristic
fixed point sets of involutions are smoothly nondisplaceable.
\end{abstract}

\section{Introduction}

Suppose that $M$ is a closed manifold and $S \subset M$ is
a subset. We say that $S$ is \emph{smoothly displaceable} if there exists a smooth time dependent vector field $X_t$ on $M$
such that its time one flow $\phi_X:=\phi^1_X$ displaces $S$ in the sense that
$$\phi_X(S) \cap S=\emptyset.$$
Otherwise we call $S$ \emph{smoothly nondisplaceable}. For example if $S$ is a closed submanifold of $M$ of half dimension
such that the homological selfintersection product of $S$ with itself does not vanish, than $S$ is smoothly nondisplaceable. One should contrast the notion of smooth displaceability with the notion of (Hamiltonian) displaceability in symplectic manifolds where the time dependent vector field is additionally required to be Hamiltonian \cite{hofer}. We recall that if $(M,\omega)$ is a symplectic manifold and $H \colon M \to \mathbb{R}$ is a smooth function referred to as the Hamiltonian than its associated Hamiltonian vector field is defined implicitly by the requirement
$$dH=\omega(\cdot,X_H).$$
We denote the time one flow of the Hamiltonian vector field of a maybe additionally time dependent Hamiltonian by
$$\phi_H:=\phi_{X_H}$$
and refer to it as a Hamiltonian diffeomorphism. It has the property that it preserves the symplectic form, i.e.,
\begin{equation}\label{sympl}
\phi_H^* \omega=\omega.
\end{equation}
Of course if a subset in a symplectic manifold is Hamiltonian displaceable it is as well smoothly displaceable. The converse is not true. The simplest example is to consider the equator on the two dimensional sphere. The equator is smoothly displaceable but it cannot be Hamiltonianly displaced because a Hamiltonian diffeomorphism is area preserving by (\ref{sympl}).
\\ \\
In this note we prove the following result.

\begin{theorem}\label{main}
Assume that the Euler characteristic of $M$ is odd and $\rho  \in \mathrm{Diff}(M)$ is a smooth involution, i.e.,
$\rho^2=\mathrm{id}|_M$, then its fixed point set
$\mathrm{Fix}(\rho) \subset M$ is smoothly nondisplaceable.
\end{theorem}

To convince ourself that the assertion of the theorem can fail if the Euler characteristic of $M$ is even, we can take again the simple example of the equator on the two dimensional sphere which is a fixed point set of a smooth involution although it is smoothly displaceable. Nevertheless it is Hamiltonianly nondisplaceable. This is an instance of the Arnold-Givental conjecture \cite{givental} which we recall next.
A \emph{real structure} on a symplectic manifold $(M,\omega)$
is an antisymplectic involution, i.e., a diffeomorphims $\rho \colon M \to M$ satisfying
$$\rho^2=\mathrm{id}|_M, \quad \rho^* \omega=-\omega.$$
The fixed point set
$$L:=\mathrm{Fix}(\rho) \subset M$$
is a (maybe empty) Lagrangian submanifold of $M$. The triple
$(M,\omega,\rho)$ is referred to as a real symplectic manifold.
We can now formulate the Arnold-Givental conjecture
\begin{conj}[Arnold-Givental conjecture]
Assume that $(M,\omega,\rho)$ is a closed real symplectic manifold and $\phi_H$ is a Hamiltonian diffeomorphism with the property that $L$ and $\phi_H(L)$ intersect transverally, then
the number of intersection points of $L$ and $\phi_H(L)$ can
be estimated from below in topological terms by the sum of $\mathbb{Z}_2$-Betti numbers of $L$, i.e.,
$$\#(L \cap \phi_H(L)) \geq \sum_{i=0}^{\mathrm{dim}(L)}b_i(L;\mathbb{Z}_2).$$
\end{conj}
If the Arnold-Givental conjecture holds true we get as an immediate Corollary the following weak version of the Arnold-Givental conjecture.
\begin{conj}[Weak Arnold-Givental conjecture]
Under the assumption that $(M,\omega,\rho)$ is a closed real symplectic manifold with the property that $L=\mathrm{Fix}(\rho) \neq \emptyset$, the Lagrangian submanifold $L$ is Hamiltonianly nondisplaceable. 
\end{conj}
As a consequence of Theorem~\ref{main} we obtain
\begin{cor}
For real symplectic manifolds of odd Euler characteristic
the weak Arnold-Givental conjecture is true.
\end{cor}
However, we point out that in the case of odd Euler characteristic the weak Arnold-Givental conjecture is by Theorem~\ref{main} a smooth and not a symplectic phenomenon. On the other hand for manifolds of even Euler characteristic the symplectic assumptions cannot be disposed of as the example of the equator on the two dimensional sphere shows. Although there are special cases of the Arnold-Givental conjecture proven, see for example \cite{frauenfelder} and the literature cited therein, the full version of the conjecture is to the knowledge of the author still open. Even the weak version does not seem to be known in general. 
\\ \\
As a further nonsymplectic Corollary of Theorem~\ref{main} we obtain the following result first proved by Conner and Floyd in 
\cite[p.\,71]{conner-floyd} using cobordism theory and later reproved by Bredon in \cite{bredon} by cohomological methods.

\begin{cor}
Assume that $M$ is a closed manifold of odd Euler characteristic and $\rho$ is a smooth involution on $M$. Then at least one of the connected components of the fixed point set
$\mathrm{Fix}(\rho)$ has dimension greater or equal half the dimension of $M$.
\end{cor}
\textbf{Acknowledgements:} The author was supported by DFG grant FR 2637/2-1.

\section{Proof of the Theorem}

We assume that $M$ is a closed manifold, $\rho \colon M \to M$
is a smooth involution, and $\phi \colon M \to M$ is a diffeomorphism. Our trick is to examine the fixed point set of
the diffeomorphism
$$\psi=\phi \rho \phi^{-1}\rho \in \mathrm{Diff}(M).$$
We use the abbreviation
$$L:=\mathrm{Fix}(\rho).$$
We need the following Lemma.
\begin{lemma}\label{free}
The fixed point set of $\psi$ is invariant under $\rho$. Moreover, if $\phi$ displaces $L$ from itself, i.e., $\phi(L) \cap L=\emptyset$, then the induced involution of $\rho$ on $\mathrm{Fix}(\psi)$ is free. 
\end{lemma}
\textbf{Proof: }We first check that $\mathrm{Fix}(\psi)$ is invariant under $\rho$.
Hence suppose that $x \in \mathrm{Fix}(\psi)$, i.e.,
$$x=\psi(x).$$
Equivalently,
\begin{equation}\label{fix}
x=\psi^{-1}(x)=\rho \phi \rho \phi^{-1}(x)
\end{equation}
where for the last equality we used the fact that $\rho^{-1}=\rho$ because $\rho$ is an involution. We compute
$$\psi(\rho(x))=\phi \rho \phi^{-1} \rho^2(x)=\phi \rho \phi^{-1}(x)=\rho^2 \phi \rho \phi^{-1}(x)
=\rho(x)$$
implying that $\rho(x)$ is a fixed point of $\psi$ as well. In particular, $\mathrm{Fix}(\psi)$ is invariant
under the involution $\rho$. 
\\ \\
We now suppose that
\begin{equation}\label{noelement}
\phi(L) \cap L=\emptyset
\end{equation}
and
$$x \in \mathrm{Fix}(\psi).$$
In order to prove the second part of the lemma our task is to show that $x \neq \rho(x)$, i.e.,
\begin{equation}\label{neithernor}
x \notin L=\mathrm{Fix}(\rho).
\end{equation}
To prove (\ref{neithernor}) we argue by contradiction and assume that 
$x \in L$, i.e.,
$$x=\rho(x).$$
Combining this with (\ref{fix}) we compute
$$x=\rho(x)=\phi \rho \phi^{-1}(x)$$
or equivalently
$$\phi^{-1}(x)=\rho \phi^{-1}(x)$$
implying that $\phi^{-1}(x) \in \mathrm{Fix}(\rho)=L$ or in other words $x \in \phi(L)$. That means
that $x$ belongs to $L$ as well as to $\phi(L)$ in contradiction to (\ref{noelement}). This shows the truth of
(\ref{neithernor}) and finishes the proof of the lemma. \hfill $\square$
\\ \\
In order to see how Lemma~\ref{free} implies Theorem~\ref{main}
we have to recall the Lefschetz-Hopf theorem, see \cite[Proposition VII.6.6]{dold}. This theorem tells us that the Lefschetz number of a continuous map coincides with the sum of the indices at the fixed points. Now suppose we have a diffeomorphism $\psi \colon M \to M$ which is homotopic to the identity. Then its Lefschetz number is just the Euler characteristic of $M$. Moreover, suppose that $\psi$ is \emph{nondegenerate} in the sense that for every fixed point
$x \in \mathrm{Fix}(\psi)$ it holds that
$$\mathrm{det}(d \psi(x)-\mathrm{id}|_{T_x M}) \neq 0,$$
i.e., $1$ is not an eigenvalue of the differential of $\psi$ at every fixed point. In this case the index of $\psi$ at each fixed point is either one or minus one. If we count modulo two we do not need to care about signs and therefore obtain the following consequence of the Lefschetz-Hopf theorem
\begin{lemma}\label{lef}
Suppose that $M$ is a closed manifold and $\psi \in \mathrm{Diff}(M)$ is a nondegenerate diffeomorphism homotopic to the identity, then
$$\#_2 \mathrm{Fix}(\psi)=\chi_2(M)$$
where $\#_2$ denotes the cardinality modulo two and $\chi_2$
denotes the Euler characteristic modulo two. 
\end{lemma}
Armed with Lemma~\ref{free} and Lemma~\ref{lef} we are now in position to prove Theorem~\ref{main}.
\\ \\
\textbf{Proof of Theorem~\ref{main}:} We argue by contradiction 
and assume that $L=\mathrm{Fix}(\rho)$ is smoothly displaceable. This means that there exists $\phi \in \mathrm{Diff}(M)$ homotopic to the identity such that 
$\phi(L)\cap L=\emptyset$. Because $\phi$ is homotopic to the identity it follows that $\psi=\phi \rho \phi^{-1} \rho$ is
homotopic to the identity as well. Now first suppose that
$\psi$ is nondegenerate. In this case it follows from 
Lemma~\ref{lef} and the assumption of the theorem that the Euler characteristic of $M$ is odd that the number of fixed points of $\psi$ is odd. This however contradicts Lemma~\ref{free} which tells us that there exists a free $\mathbb{Z}_2$-action on $\mathrm{Fix}(\psi)$. This contradiction proves the theorem in the case that $\psi$ is nondegenerate.

For the general case we note that because $L$ is a closed submanifold of $M$, the property of displacing $L$ is open in the space of diffeomorphisms of $M$. Because for a generic diffeomorphism $\phi$ the corresponding diffeomorphism
$\psi$ is nondegenerate the general case follows from the special case by maybe slightly perturbing $\phi$. This finishes the proof of the theorem. \hfill $\square$

\end{document}